# STOCHASTIC SPATIAL MODELS OF HOST-PATHOGEN AND HOST-MUTUALIST INTERACTIONS I

BY N. LANCHIER AND C. NEUHAUSER[1]

*CNRS–Université de Rouen and University of Minnesota*

Mutualists and pathogens, collectively called symbionts, are ubiquitous in plant communities. While some symbionts are highly host-specific, others associate with multiple hosts. The outcomes of multispecies host-symbiont interactions with different degrees of specificity are difficult to predict at this point due to a lack of a general conceptual framework. Complicating our predictive power is the fact that plant populations are spatially explicit, and we know from past research that explicit space can profoundly alter plant-plant interactions. We introduce a spatially explicit, stochastic model to investigate the role of explicit space and host-specificity in multispecies host-symbiont interactions. We find that in our model, pathogens can significantly alter the spatial structure of plant communities, promoting coexistence, whereas mutualists appear to have only a limited effect. Effects are more pronounced the more host-specific symbionts are.

**1. Introduction.** The diversity and structure of plant communities are largely determined by nutrient availability, competition among plants, herbivory and associations between plants and their symbionts. The first three have been the focus of much ecological research, both empirical and theoretical. The role of symbionts on diversity and structure has received less attention [1] and will be the topic of this paper.

A symbiont is an organism that lives in close association with a host. It can have either a harmful effect, in which case we call it a pathogen, or a beneficial effect, in which case we call it a mutualist. Symbionts are ubiquitous. For instance, more than 90% of terrestrial plants associate with mycorrhizal fungi [18], a beneficial association that supplies nutrients to the plant and, in

Received April 2005; revised August 2005.
[1]Supported in part by NSF Grants DMS-00-72262 and DMS-00-83468.

return, carbon to the fungal partner. Fecundity and viability of virtually all plants are affected by pathogens, sometimes with devastating effects, such as Dutch elm disease or chestnut blight that point to the important role of pathogens in structuring plant communities. Mathematical models play an important role in elucidating the roles of symbionts in community dynamics.

Modeling of disease dynamics has had a long tradition, starting with the model by Kermack and McKendrick [14], which describes the course of a disease outbreak caused by a single disease infecting a single host. This model and its extensions have yielded enormously valuable insights into disease dynamics and potential control strategies. Although originally developed to describe epidemics in human populations, it can equally well be applied to plant diseases. A key concept of disease dynamics is the basic reproductive rate $R_0$, which is defined as the expected number of secondary infections caused by an infected individual when introduced into a population of susceptible individuals [2]. The condition for a disease outbreak is given in the biological literature as $R_0 > 1$. This condition is based on a single-host, single-disease model in a nonspatial population. Including spatial structure in the form of local interactions has shown that for a disease to spread, $R_0$ needs to exceed a threshold that is greater than that for a nonspatial population. The reason for this is the lack of sufficient numbers of susceptible individuals near the location of a disease outbreak once the disease starts spreading. One of the first models where this has been demonstrated mathematically is the contact process [11, 16].

Much of the theoretical work in the epidemiological literature has focused on single-host, single-disease dynamics. A rapidly increasing empirical body of work on multispecies host-disease dynamics necessitates the development of a theoretical framework. This has only begun recently (see [13], and references therein). An attempt for a broad classification was made by Holt et al. [13] using a graphical isocline framework that allows for generalizations that are applicable to a wide range of host-pathogen models. A recent model by Dobson [4] investigates persistence of a pathogen that can infect multiple hosts.

Modeling of host-mutualist interactions has primarily focused on evolutionary questions, such as the evolution of cheaters (i.e., symbionts that receive benefits but do not confer them). Almost no modeling has been done on the effects of local (spatial) interactions on host-mutualist dynamics. Similarly, no theoretical framework has been developed for multispecies host-mutualist interactions.

An important component of multihost, multisymbiont models is the degree to which different symbionts and hosts can associate with each other. This is referred to as *specificity*. A *specialist* symbiont associates with a very small number of hosts; a *generalist* symbiont associates with many hosts.

The ease of transmission of a symbiont to a host, referred to as *transmissibility*, is another important factor in host-symbiont interactions.

In this paper we investigate the role of spatial structure caused by local interactions, such as symbiont transmission and host dispersal, on persistence of host-symbiont associations for both generalists and specialists in multihost, multisymbiont models. We employ the simplest of all multispecies models to describe the host dynamics, the voter model [3, 12]. The voter model is defined on the $d$-dimensional integer lattice, where each lattice site is occupied by an individual characterized by one of a finite number of types. Individuals give birth to offspring of their own kind at a constant rate, and their offspring displace randomly chosen individuals within their dispersal neighborhood. The dynamics imply that all sites remain occupied at all times. Into this population, we introduce symbionts with varying degrees of specificity and transmissibility. More precisely, our spatial model is a continuous-time Markov process $\xi_t : \mathbb{Z}^d \longrightarrow \{1, 2, \ldots, \kappa\} \times \{0, 1, \ldots, \kappa\}$ where the integer $\kappa$ denotes both the number of hosts and the number of symbionts involved in the interaction. A site $x \in \mathbb{Z}^d$ is said to be occupied by an *unassociated* host of type $i$, $i = 1, 2, \ldots, \kappa$, if $\xi(x) = (i, 0)$, and by a host of type $i$, $i = 1, 2, \ldots, \kappa$, *associated* with a symbiont of type $j$, $j = 1, 2, \ldots, \kappa$, if $\xi(x) = (i, j)$. We will write $\xi_t(x) = (\xi_t^1(x), \xi_t^2(x))$, where $\xi_t^1(x)$ denotes the type of the host present at $x$ at time $t$ and $\xi_t^2(x)$ denotes the type of the symbiont present at $x$ at time $t$, with $\xi_t^2(x) = 0$ denoting the absence of a symbiont. We set $\|x\| = \sup_{i=1,2,\ldots,d} |x_i|$. The evolution at site $x \in \mathbb{Z}^d$ is described by the transition rates

$$(i,j) \to (k,0) \quad \text{at rate} \quad \lambda \sum_{0 < \|x-z\| \leq R_1} \left\{ \mathbb{1}_{\{\xi(z)=(k,0)\}} + g \sum_{\ell=1}^{\kappa} \mathbb{1}_{\{\xi(z)=(k,\ell)\}} \right\}$$

$$(i,0) \to (i,j) \quad \text{at rate} \quad c_{ij} \sum_{0 < \|x-z\| \leq R_2} \sum_{\ell=1}^{\kappa} \mathbb{1}_{\{\xi(z)=(\ell,j)\}}.$$

The transition $(i, j) \to (k, 0)$ is the birth of an unassociated host at $x$ by either unassociated or associated neighboring hosts. The birth rate of unassociated hosts is equal to $\lambda$. The parameter $g$ indicates the variation of the birth rate of hosts associated with a symbiont. If $0 \leq g < 1$, the symbiont is a pathogen; if $g = 1$, the symbiont has no effect on the birth rate of the host and we refer to this as the neutral case; if $g > 1$, the symbiont is a mutualist. The transition $(i, 0) \to (i, j)$ is the transmission of a neighboring symbiont $j$ to an unassociated host of type $i$ at $x$. The parameters $c_{ij}$ denote the rate at which symbiont $j$ infects host $i$. This parameter will allow us to mimic specialist and generalist symbionts. Births and infections occur within a local neighborhood, with $R_1$ denoting the birth range of hosts, and $R_2$ the infection range of symbionts. Neighborhoods are punctured boxes with side $2R_i + 1$, $i = 1, 2$, centered at site $x$, that is, $\mathcal{N}_x^i = \{z \in \mathbb{Z}^d : 0 < \|x - z\| \le R_i\}$. The cardinality of this set is denoted by $\nu_{R_i} = |\mathcal{N}_x^i|$.

Before we describe the behavior of the spatially explicit, stochastic model, we will look at the mean-field model [9]. The mean-field model is described by a system of differential equations for the densities of unassociated and associated hosts. To define it, we let $u_i$ denote the density of unassociated hosts of type $i$, $i = 1, 2, \ldots, \kappa$, and $v_{ij}$ denote the density of host $i$ associated with symbiont $j$, $i, j = 1, 2, \ldots, \kappa$. It follows from the dynamics of the spatially explicit, stochastic model that at all times

$$\sum_i u_i + \sum_{i,j} v_{ij} = 1.$$

Furthermore, we assume that for $i = 1, 2, \ldots, \kappa$, $c_{ii} = \beta$, and for $i, j = 1, 2, \ldots, \kappa$ with $i \ne j$, $c_{ij} = \alpha$ with $0 \le \alpha \le \beta$. We define

$$u_\cdot = \sum_{i=1}^\kappa u_i, \quad v_{\cdot j} = \sum_{i=1}^\kappa v_{ij}, \quad v_{i\cdot} = \sum_{j=1}^\kappa v_{ij} \quad \text{and} \quad v_{\cdot\cdot} = \sum_{i=1}^\kappa \sum_{j=1}^\kappa v_{ij}.$$

One way to obtain the mean-field limit is to set the neighborhood ranges $R_1$ and $R_2$ equal to $R$ and then let $R$ go to infinity. To obtain a meaningful limit, we also need to rescale the parameters $\lambda$, $\alpha$ and $\beta$ by the neighborhood size $\nu_R$; that is, we set $\lambda = \frac{1}{\nu_R}$ (this also sets the time scale), and define

$$\alpha = \frac{a}{\nu_R} \quad \text{and} \quad \beta = \frac{b}{\nu_R}.$$

In the limit $R \to \infty$, sites become independent. If, in addition, the initial configuration is translation invariant, the dynamics of the densities for $i \ne j$ is then described by the following system of differential equations, called mean-field equations:

$$\frac{du_i}{dt} = (1 - u_i)(u_i + gv_{i\cdot}) - u_i \sum_{j \ne i}(u_j + gv_{j\cdot}) - bu_i v_{\cdot i} - a \sum_{j \ne i} u_i v_{\cdot j},$$

$$\frac{dv_{ii}}{dt} = bu_i v_{\cdot i} - v_{ii}(u_\cdot + gv_{\cdot\cdot}),$$

$$\frac{dv_{ij}}{dt} = au_i v_{\cdot j} - v_{ij}(u_\cdot + gv_{\cdot\cdot}).$$

When $a = 0$, the symbionts are specialists. As $a$ increases to $b$, the association turns into a generalist relationship. The following results are proved in Section 2. When $g = 1$, the system has a conserved quantity, namely the initial host densities $h_i = u_i + v_{i\cdot}$, $i = 1, 2, \ldots, \kappa$. If $(\kappa - 1)a + b > \kappa$, then for $g \ne 1$, there exists a nontrivial equilibrium with $u_1 = u_2 = \cdots = u_\kappa \ge 0$ and $v_{1\cdot} = v_{2\cdot} = \cdots = v_{\kappa\cdot} > 0$ such that for $i = 1, 2, \ldots, \kappa$,

$$u_i = \frac{g}{(\kappa - 1)a + b - \kappa(1 - g)} \quad \text{and} \quad h_i = \frac{1}{\kappa}.$$

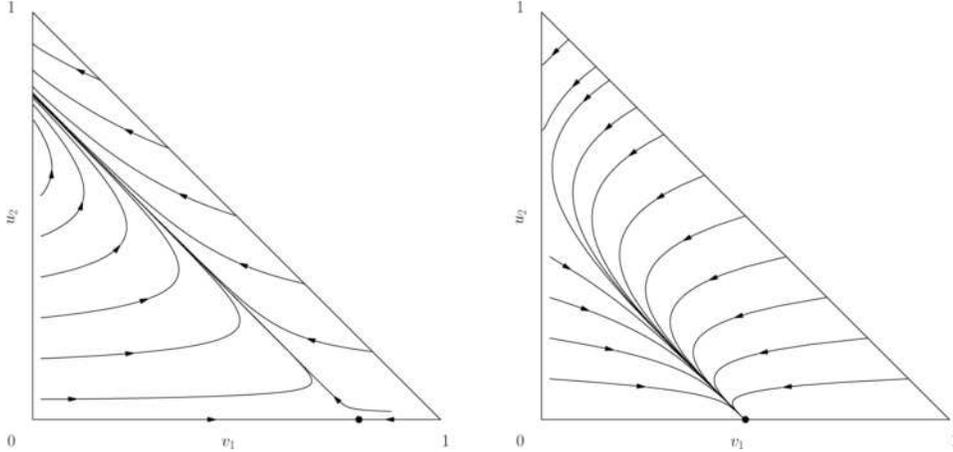

FIG. 1. *Solution curves of the mean-field model with $a = 0$ and $b = 3$. The parameter $g$ is equal to $0.5$ and $2$, respectively.*

Furthermore, for $i \neq j$

$$\frac{v_{ij}}{v_{ii}} = \frac{a}{b}.$$

Numerical simulations indicate that the nontrivial equilibrium is locally stable for $a < b$ when $g < 1$, but not for $g > 1$. In addition, if $g = 0$ and $(\kappa - 1)a + b > \kappa$, all hosts will be associated at equilibrium. If initially only two hosts and one symbiont, say symbiont 1, are present, then for $0 < g < 1$, simulations indicate that pathogen 1 will go extinct and both host 1 and host 2 may coexist. We will conjecture similar behavior for the spatial model, except in the one-dimensional, nearest-neighbor case (see Theorem 3 and discussion following the theorem). When $g > 1$, host 2 goes extinct provided the initial density of mutualists of type 1 is sufficiently large. The spatial analogue of this result is proved in Theorem 2. Both results are illustrated in Figure 1.

We now return to the spatially explicit model with parameters $\lambda$, $\alpha$ and $\beta$. To define the time scale, we set $\lambda = 1$. We will discuss both the generalist case $\alpha = \beta$ and the specialist case $\alpha = 0$, and provide comparisons with the mean-field model.

*Generalist interactions.* We consider the generalist case $\alpha = \beta$ of the spatially explicit, stochastic model. In the corresponding mean-field model, coexistence of hosts and symbionts is possible when $b > 1$. For $b \leq 1$, symbionts are unable to persist. When $a = b$, symbionts are no longer host-specific, and the mean-field model can be reduced to one with $u_i$ and $v_i$.. It is not hard to see then that the nontrivial equilibrium of unassociated hosts, $(u., v..)$, is locally stable since in this case

$$\frac{dv..}{dt} = v..[(b-1)u. - gv..].$$

With $u. + v.. = 1$, it follows that the boundary point $v.. = 0$ is unstable for $b > 1$ and that

$$(u., v..) = \left(\frac{g}{b-1+g}, \frac{b-1}{b-1+g}\right)$$

is locally stable. Furthermore, any vector $(h_1, h_2, \ldots, h_\kappa)$ with $h_i \geq 0$ and $\sum_i h_i = 1$ gives rise to an equilibrium if we set

$$u_i = \frac{g}{b-1+g} h_i \quad \text{and} \quad v_i. = \frac{b-1}{g} u_i.$$

The behavior of the spatially explicit model is more complicated and may depend on the spatial dimension, as we will see in the following.

If $g = 1$, then the symbionts have no effect on the hosts, which means that the spatially explicit processes $\xi_t^1$ and $\xi_t^2$ are stochastically independent. Moreover, by looking at the transition rates, it is easy to see that $\xi_t^1$ is a multitype voter model run at rate 1, and that $\xi_t^2$ is a multitype contact process in which particles give birth at rate $\beta \nu_{R_2}$ and die at rate $\nu_{R_1}$. See, respectively, [12, 17] for a study of these two processes. It follows that there exists a critical value $\beta_c(R_1, R_2) \in (0, \infty)$ that depends on the neighborhood sizes $\nu_{R_1}$ and $\nu_{R_2}$ such that the symbionts can survive if and only if $\beta > \beta_c(R_1, R_2)$. If we ignore host and symbiont types but rather focus on associated versus unassociated hosts, then for $\beta > \beta_c(R_1, R_2)$, regardless of the spatial dimension, there exists a nontrivial stationary distribution of associated and unassociated hosts. Moreover, if $d \geq 3$, there exists a stationary distribution in which all hosts and symbionts coexist.

Unfortunately, we cannot say much about coexistence when $g \neq 1$. To analyze this case, we define the "color-blind" process where a site is in state 0 if it is occupied by an unassociated host, and in state 1 if it is occupied by an associated host. We obtain a particle system $\zeta_t : \mathbb{Z}^d \longrightarrow \{0, 1\}$ with transitions at $x \in \mathbb{Z}^d$:

$$0 \to 1 \quad \text{at rate} \quad \beta \sum_{0 < \|x-z\| \leq R_2} \mathbb{1}_{\{\zeta(z)=1\}},$$

$$1 \to 0 \quad \text{at rate} \quad \sum_{0 < \|x-z\| \leq R_1} \{\mathbb{1}_{\{\zeta(z)=0\}} + g \mathbb{1}_{\{\zeta(z)=1\}}\}.$$

When $g = 0$, the process reduces to a biased voter model. When $g = 1$, it reduces to a contact process with birth rate $\beta \nu_{R_2}$ and death rate $\nu_{R_1}$. We denote the critical value of this contact process by $\beta_c(R_1, R_2)$ as above. A standard coupling argument allows us to compare the processes with $g \neq 1$ and $g = 1$, and to deduce that if $g \leq 1$ and $\beta > \beta_c(R_1, R_2)$, then $\zeta_t$ has a nontrivial stationary measure, while if $g \geq 1$ and $\beta \leq \beta_c(R_1, R_2)$, then the mutualists die out, that is, $\zeta_t \Rightarrow \delta_0$, the "all 0" configuration [i.e., $\zeta(x) \equiv 0$]. To cover the remaining cases, we introduce the contact process $\eta_t$ in which particles give birth at rate $\beta \nu_{R_2}$ and die at rate $g \nu_{R_1}$. Then $\eta_t$ has a nontrivial stationary measure if and only if $\beta > g\beta_c(R_1, R_2)$ which, with a new coupling argument, implies that if $g \leq 1$ and $\beta \leq g\beta_c(R_1, R_2)$, then the pathogens die out, while if $g \geq 1$ and $\beta > g\beta_c(R_1, R_2)$, then $\zeta_t$ has a nontrivial stationary measure.

We now focus on the case $\beta > 1$ and $g > 0$ close to 0. First of all, we observe that if $g = 0$ and $R_1 = R_2$, then the process $\zeta_t$ is the biased voter model with parameters $\beta$ and 1. In particular, $P(\zeta_t(x) = 0) \to 1$ if $\beta < 1$ while $P(\zeta_t(x) = 1) \to 1$ if $\beta > 1$ provided we start with infinitely many 0's and 1's at time 0. Moreover, in the latter case, fixation occurs for the process $\xi_t$ since hosts associated with pathogens are now sterile. The behavior is identical to that of the mean-field model. We will use a perturbation argument in Section 4 to show that if $\beta > 1$ and $g > 0$ is sufficiently close to 0, then the pathogens still survive. The results are summarized in Figure 2 and in the following theorem where "$\Rightarrow$" denotes weak convergence and $\delta_0$ is the distribution that concentrates on the "all 0" configuration.

THEOREM 1. *Assume that $\alpha = \beta$ and that $\zeta_0$ is translation invariant with $P(\zeta_0(x) = 1) > 0$.*

(a) *If $g \leq 1$, then $\zeta_t \Rightarrow \delta_0$ if $\beta \leq g\beta_c$, and a nontrivial equilibrium exists if $\beta > \beta_c$. If $g \geq 1$, then $\zeta_t \Rightarrow \delta_0$ if $\beta \leq \beta_c$, and a nontrivial equilibrium exists if $\beta > g\beta_c$.*

(b) *If $\beta > 1$, there exists $g_c > 0$ such that if $g \leq g_c$, then $\zeta_t \Rightarrow \mu$ with $\mu(\zeta(x) = 1) \neq 0$.*

Part (b) of this theorem will be proved in Section 4.

*Specialist interactions.* In the specialist case $\alpha = 0$ and $\beta > 0$, the process is more difficult to investigate since the evolution of each symbiont strongly depends on the configuration of the host population. That is, there is no particle system $\zeta_t : \mathbb{Z}^d \longrightarrow \{0, 1\}$ which allows us to describe the global evolution of the symbionts regardless of their type. Since for any $i = 1, 2, \ldots, \kappa$ the symbiont $i$ can live only through hosts of type $i$, it is, however, easy to deduce from a coupling argument that if $g \leq 1$, then the processes with $\alpha = 0$ and $\alpha = \beta$ can be defined on the same space so that, starting from the same configuration, the process with $\alpha = 0$ has fewer pathogens. In words, the survival of the pathogens is harder to obtain with specialist interactions. In particular, if $g \leq 1$ and $\beta \leq \max(\beta_c g, 1)$, then the pathogens die out.

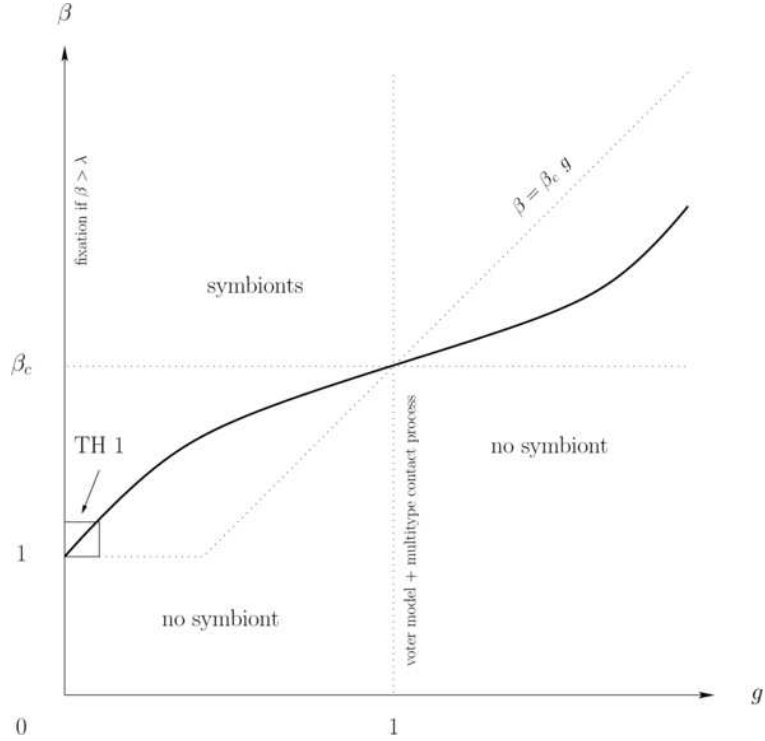

Fig. 2. *Phase diagram of the spatial model with generalist interactions.*

The next step is to extend the results of the mean-field model summarized in Figure 1 to the corresponding spatial model with short-range interactions. To do this, we consider a population of two hosts with only one type of symbiont, say symbiont of type 1, and start the evolution with all the hosts of type 1 associated with a symbiont. Then, in the limiting case $\beta = \infty$, the transition $(1,0) \to (1,1)$ is instantaneous; that is, unassociated hosts of type 1 become instantaneously associated with a symbiont, provided that $R_1 \leq R_2$ to avoid the problem of an *isolated* unassociated host that cannot be reached by any symbiont. This implies that the process $\xi_t^1$ is a biased voter model in which hosts of type 1 give birth at rate $g\nu_{R_1}$ and hosts of type 2 at rate $\nu_{R_1}$. In particular, if $g < 1$, then $\xi_t^1 \Rightarrow \delta_2$, the "all 2" configuration, while if $g > 1$, then $\xi_t^1 \Rightarrow \delta_1$, the "all 1" configuration. Theorem 2 tells us that, in any dimension, the result still holds if $g > 1$ and $\beta \in (0, \infty)$ is sufficiently large.

THEOREM 2. *Assume that $\alpha = 0$, $g > 1$ and $R_1 \leq R_2$. At time 0, $\xi_0(x) = (1,1)$ or $(2,0)$ for all $x \in \mathbb{Z}^d$. Then, there is a critical value $\beta_{cr}^{\mathrm{Th}2} \in (0, \infty)$ such that if $\beta > \beta_{cr}^{\mathrm{Th}2}$, then $\xi_t^1 \Rightarrow \delta_1$.*

The arguments in our proof, however, do not imply the analogous result for $g < 1$. We think actually that except in the one-dimensional nearest-neighbor case, $P(\xi_t^1(x) = 2) \not\to 1$. Before explaining our intuition, we describe the behavior of the one-dimensional process with nearest-neighbor interactions for the hosts and short-range interactions for the pathogens. Assume that at time 0, $\xi_0^1(x) = 1$ for $x > 0$ with infinitely many pathogens of type 1, and $\xi_0^1(x) = 2$ for $x \leq 0$ with no associated pathogens. Denote by $r_t^2 = \sup\{x \in \mathbb{Z} : \xi_t^1(x) = 2\}$ the rightmost host of type 2. Then $r_0^2 = 0$. The following result implies that for $\beta$ large enough, 2's will spread to the right and eliminate 1's together with their associated pathogens.

THEOREM 3. *Assume $d = 1$ and $R_1 = R_2 = 1$. If $\alpha = 0$ and $g < 1$, there exists $\beta_{cr}^{\mathrm{Th}3} \in (0, \infty)$ such that if $\beta > \beta_{cr}^{\mathrm{Th}3}$, then $r_t^2 \to \infty$ as $t \to \infty$ with probability* 1.

We conjecture that this result should only be true in the one-dimensional nearest-neighbor case. Here is the intuition. We first observe that except in the one-dimensional nearest-neighbor case, the dynamics produce *isolated* hosts, that is, hosts of type 1 (resp. 2) surrounded by a cluster of hosts of type 2 (resp. 1). As suggested by Theorem 2, when $g > 1$, isolated 2's are *swallowed* very quickly by surrounding 1's. On the other hand, when $g < 1$, an isolated 1 located in a linearly growing cluster of 2's cannot be invaded anymore by any pathogen as soon as the cluster has reached some critical size. In this expanding region, the process then behaves like an unbiased voter model in which 1's can now compete with 2's, and survive. See Figure 3 in Section 5 for simulations in $d = 2$.

Finally, we investigate the coexistence of symbionts in the neutral case $g = 1$. We observe that, in this case, the first coordinate process $\xi_t^1$ performs a voter model run at rate 1. In particular, in $d \leq 2$, there does not exist any stationary distribution $\mu$ such that $\mu(\xi^1(x) = i) \neq 0$ for any $i \in \{1, 2, \ldots, \kappa\}$ (see [12]). Since in the specialist case, $\xi_t^2(x) = i$ implies that $\xi_t^1(x) = i$, the same conclusion holds for the symbionts. In $d \geq 3$, coexistence occurs for the process $\xi_t^1$, that is, there is a stationary measure $\mu$ which satisfies the condition above. However, due to the formation of clusters, the problem of coexistence of the symbionts remains a difficult question. Namely, the voter model $\xi_t^1$ exhibits clusters whose diameter can exceed some critical size, which prevents the symbionts from spreading out. To get around this difficulty, we introduce a modification of the particle system, denoted by $\hat{\xi}_t$, in which the symbionts evolve as previously but where the hosts now perform a threshold $\theta$ voter model. More precisely, the process $\hat{\xi}_t$ evolves according to the following transitions at $x \in \mathbb{Z}^d$:

$$(i, j) \to (k, 0) \quad \text{at rate} \quad \begin{cases} 1, & \text{if card}\{z \in \mathbb{Z}^d : 0 < \|x - z\| \leq R_1 \\ & \qquad \text{and } \hat{\xi}^1(z) = k\} \geq \theta, \\ 0, & \text{otherwise,} \end{cases}$$

$(i, 0) \to (i, i)$  at rate  $\beta \operatorname{card}\{z \in \mathbb{Z}^d : 0 < \|x - z\| \leq R_2 \text{ and } \hat{\xi}^2(z) = i\}.$

The introduction of this particle system is motivated by Theorem 1 of [7] which implies that the threshold $\theta$ voter model has a nontrivial stationary distribution which is close enough to a product measure to produce our next result. We will prove the following result.

THEOREM 4. *Let $g = 1$ and $\theta < \nu_{R_1}/\kappa$. If $R_1$ and $R_2$ are sufficiently large, then there is a critical value $\beta_{cr}^{\mathrm{Th}4} \in (0, \infty)$, depending on $\kappa$, such that the following holds: If $\beta > \beta_{cr}^{\mathrm{Th}4}$, then coexistence occurs, and if $\beta < \beta_{cr}^{\mathrm{Th}4}$, then coexistence is not possible.*

Here coexistence means that there is a stationary measure $\mu$ such that $\mu(\hat{\xi}^2(x) = i) \neq 0$ for any type $i \in \{1, 2, \ldots, \kappa\}$. Unfortunately, we do not know how to prove something better than *coexistence is not possible* when $\beta < \beta_{cr}^{\mathrm{Th}4}$. However, we conjecture that starting from product measure in which each host is associated with a symbiont and has density $1/\kappa$, all the symbionts die out. To justify our conjecture, we observe that the processes $\eta_t^i$ defined by $\eta_t^i(x) = 1$ if $\hat{\xi}_t^2(x) = i$ and $\eta_t^i(x) = 0$ otherwise, do not interact since they are confined to their associated hosts. Since Theorem 1 of [7] tells us that the hosts coexist with density of each type close to $1/\kappa$, each symbiont should remain subcritical. Finally, since the symbionts can only spread out through their host, we conjecture that $\beta_{cr}^{\mathrm{Th}4}$ is increasing with respect to $\kappa$.

*Comparison of the spatially explicit and the mean-field model.* Numerical simulations of the mean-field model indicate that coexistence is only possible when $g < 1$. Simulations of the spatially explicit model show similar behavior. When $g < 1$ and $\alpha < \beta$, then coexistence of hosts and pathogens is possible. We observed that in this case, cluster size is limited by the presence of pathogens: In the absence of pathogens, clusters grow at the expense of neighboring clusters that contain symbionts. Upon invasion by the preferred symbionts (those with infection rate $\beta$), the clusters appear to shrink again. The case $g > 1$ and $\alpha < \beta$ is quite different. Clusters of hosts with their preferred mutualists form and appear to continue to grow, just as in the voter model case. Less preferred mutualists (those with infection rate $\alpha$) do not seem to be able to persist with preferred mutualists, just as is the mean-field case of one host and two symbionts with infection rates $a$ and $b$, respectively ($a < b$). In summary, pathogens have the ability to alter the spatial structure of their hosts by promoting local diversity, whereas mutualists do not alter the spatial structure of their hosts. This difference in behavior is more pronounced the more host-specific the symbionts are.

The rest of this paper is devoted to proofs. In Section 2 we will investigate the mean-field model. In Section 3 we will prove a preliminary result about the biased voter model to prepare the proofs of Theorems 1 and 2 which will be carried out in Sections 4 and 5, respectively. Section 6 will be devoted to the proof of Theorem 3. Finally, we will investigate the coexistence of symbionts and prove Theorem 4 in Section 7.

**2. The mean-field model.** The mean-field model was introduced in Section 1. Our first claim was that the host density $h_i = u_i + v_i.$ is a conserved quantity when $g = 1$. A straightforward calculation shows that if $g = 1$,

$$\frac{d}{dt}(u_i + v_i.) = 0,$$

from which our claim follows.

We summarize the behavior of the mean-field model in the following proposition.

PROPOSITION 2.1. *For $(\kappa - 1)a + b > \kappa$ and $g \geq 0$, there exists a nontrivial equilibrium with $u_1 = u_2 = \cdots = u_\kappa \geq 0$ and $v_{1.} = v_{2.} = \cdots = v_{\kappa.} > 0$ such that for $i = 1, 2, \ldots, \kappa$,*

$$u_i = \frac{g}{(\kappa - 1)a + b - \kappa(1 - g)} \quad and \quad h_i = \frac{1}{\kappa}.$$

*Furthermore, for $i \neq j$*

$$\frac{v_{ij}}{v_{ii}} = \frac{a}{b}.$$

PROOF. If we denote by $h_i = u_i + v_i.$ the density of host $i$ (both associated and unassociated), then

$$\frac{dh_i}{dt} = u_i + gv_i. - h_i(u. + gv..).$$

By setting the right-hand side equal to 0, we obtain

$$h_i = \frac{u_i + gv_i.}{u. + gv..}.$$

It follows that

$$\frac{h_i}{h_j} = \frac{u_i(1 - g) + gh_i}{u_j(1 - g) + gh_j}$$

from which we conclude that

$$\frac{h_i}{h_j} = \frac{u_i}{u_j} = \frac{v_i.}{v_j.}.$$

In the symmetric case, $h_1 = h_2 = \cdots = h_\kappa = \frac{1}{\kappa}$, we find $u_1 = u_2 = \cdots = u_\kappa$ and $v_{1\cdot} = v_{2\cdot} = \cdots = v_{\kappa\cdot}$. The nontrivial equilibrium can then be computed explicitly. We find

$$u_i = \frac{g}{(\kappa - 1)a + b - \kappa(1 - g)} \quad \text{and} \quad h_i = \frac{1}{\kappa}.$$

Specifically, when $g = 0$, $u_i = 0$ and consequently all hosts will be associated at equilibrium. The condition for the existence of a nontrivial point equilibrium, namely $(\kappa - 1)a + b > \kappa$, follows directly from requiring that $u_i < 1/\kappa$ and $h_i = 1/\kappa$. Furthermore, it follows from

$$bu_i v_{\cdot i} = v_{ii}(u_\cdot + gv_{\cdot\cdot}) \quad \text{and} \quad au_i v_{\cdot j} = v_{ij}(u_\cdot + gv_{\cdot\cdot})$$

that

$$\frac{a}{b} = \frac{v_{\cdot i}}{v_{\cdot j}} \frac{v_{ij}}{v_{ii}}.$$

Since $v_{\cdot i} = v_{\cdot j}$ by symmetry, the last claim follows as well. $\square$

**3. Preliminary results about the biased voter model.** As explained in the Introduction, if $\alpha = \beta$ and $R_1 = R_2$, then the "color-blind" process $\zeta_t$ performs a biased voter model when $g = 0$. If we set $\alpha = 0$ and consider a population of two host types with only one symbiont type, then the process $\xi_t^1$, which describes the evolution of both host types, performs a biased voter model in the limiting case $\beta = \infty$. So, to prove Theorems 1 and 2, we will start by proving a general result about the biased voter model, and then apply a perturbation argument to extend this result to the region $g > 0$ small in the first case, and to the region $\beta < \infty$ large in the second case. Let $\beta_1, \beta_2 \in (0, \infty)$, and $\eta_t : \mathbb{Z}^d \longrightarrow \{1, 2\}$ be the biased voter model with parameters $\beta_1$ and $\beta_2$, that is, the process whose state at site $x$ changes as follows:

$$i \to j \quad \text{at rate} \quad \beta_j \sum_{0 < \|x - z\| \le R_1} \mathbb{1}_{\{\eta(z) = j\}}.$$

It is a well-known fact that if $\beta_1 > \beta_2$, then $P(\eta_t(x) = 1) \to 1$ as $t \to \infty$, provided that at time 0, the process has infinitely many 1's and 2's (see, e.g., [6], Chapter 3).

To prove Theorems 1 and 2, we will follow the strategy described in [10], Section 3. We begin with a rescaling argument to estimate the rate of convergence of $P(\eta_t(x) = 1)$. This estimate will have to be good enough so that a perturbation argument can be applied. The basic idea is to show that for given $\varepsilon > 0$, members of the family of processes under consideration, when viewed on suitable length and time scales, dominate an $M$-dependent oriented percolation process in which sites are open with probability $1 -$ $\varepsilon$ ([8], Section 4). To compare the process with a percolation process, we consider a positive integer $L$ to be fixed later, and scale space by setting

$$B = [-L, L]^d, \qquad \Phi(z) = Lz, \qquad B_z = \Phi(z) + B$$

for any $z \in \mathbb{Z}^d$. Let $\Gamma$ be a positive integer, and say that $(z, n)$ is *occupied* if all sites in $B_z$ are occupied by 1's at time $n\Gamma L$. The first step in proving Theorems 1 and 2 is the following.

PROPOSITION 3.1. *Let $\varepsilon > 0$ and $\beta_1 > \beta_2$. Then $M$, $L$ and $\Gamma$ can be chosen in such a way that the set of occupied sites dominates the set of open sites in an $M$-dependent oriented site percolation process where sites are open with probability $p = 1 - 2\varepsilon/3$.*

The key to the proof is duality ([6], Chapter 3). To define the dual process of the biased voter model, we consider two collections of independent Poisson processes $\{T_n^{x,z} : n \geq 1\}$ and $\{U_n^{x,z} : n \geq 1\}$ where $0 < \|x - z\| \leq R_1$, with parameter $\beta_2$ and $\beta_1 - \beta_2$, respectively. At times $T_n^{x,z}$ we draw an arrow from $z$ to $x$ and put a $\delta$ at site $x$, while at times $U_n^{x,z}$ we draw an arrow from $z$ to $x$ without putting a $\delta$ at the tip. The process is then obtained from the graphical representation as follows: At time $T_n^{x,z}$, the particle at $x$ imitates the one at $z$. At time $U_n^{x,z}$, the site $x$ becomes occupied by a particle of type 1 if $z$ is. We say that there is a *path* from $(x, 0)$ to $(z, t)$ if there is a sequence of times $s_0 = 0 < s_1 < \cdots < s_{n+1} = t$ and spatial locations $x_0 = x, x_1, \ldots, x_n = z$ such that the following two conditions hold:

1. For $i = 1, 2, \ldots, n$, there is an arrow from $x_{i-1}$ to $x_i$ at time $s_i$.
2. For $i = 0, 1, \ldots, n$, the vertical segments $\{x_i\} \times (s_i, s_{i+1})$ do not contain any $\delta$'s.

Finally, we say that there exists a *dual path* from $(x, t)$ to $(z, t-s)$, $0 \leq s \leq t$, if there is a path from $(z, t-s)$ to $(x, t)$, and define the *dual process starting at* $(x, t)$ by setting

$$\hat{\eta}_s^{(x,t)} = \{z \in \mathbb{Z}^d : \text{there is a dual path from } (x, t) \text{ to } (z, t-s)\}$$

for any $0 \leq s \leq t$. The reason why we introduce the dual process is that it allows us to deduce the state of site $x$ at time $t$ from the configuration at earlier times. More precisely,

$$\eta_t(x) = 1 \quad \text{if and only if} \quad \eta_{t-s}(z) = 1 \quad \text{for some } z \in \hat{\eta}_s^{(x,t)}.$$

See [6], Chapter 3. The strategy to proving Proposition 3.1 can then be summarized as follows: Let $T = \Gamma L$ and $x \in B_z$ with $\|z\| = 1$. Then, we will prove that, with probability arbitrarily close to 1, there exists a dual path $A_s$ starting at $(x, T)$ and landing in the target set $B$. More precisely, we will prove that $A_s$ hits the set $J = [-R_1, R_1]^d$ by time $T$ where $R_1 < L/2$, and then stays inside $B$ until time $T$. Recall that $R_1$ denotes the range of the interactions. In particular, if $B$ is void of 2's at time 0, then, with probability close to 1, $B_z$ will be void of 2's as well $T$ units of time later. To define the dual path $A_s$, we start the process at $A_0 = (x, T)$ and go down the graphical representation. If $A_s$ comes across a $\delta$ at some time $s = T - T_n^{x,z}$ with $x = A_s$, then move $A_s$ to $z$. If $A_s$ meets the tip of an arrow that is without a $\delta$ at some time $s = T - U_n^{x,z}$, then move $A_s$ to $z$ only if it takes it closer to 0. Intuitively, this should cause $A_s$ to drift toward the set $B$. We now make this argument precise in a series of lemmas.

LEMMA 3.2. *Assume that $x \in B_z$, $\|z\| = 1$ and $\beta_1 > \beta_2$. There exist $C_1$, $\gamma_1 \in (0, \infty)$ such that*

$$\sup_{x \in B_z} P_x(A_s \notin J \text{ for all } s \leq T) \leq C_1 \exp(-\gamma_1 L)$$

*for $L$ and $\Gamma$ sufficiently large. Here, the subscript $x$ indicates the starting point.*

PROOF. Let $\sigma_k$ denote the $k$th time $A_s$ encounters the tip of an arrow (with or without a $\delta$). At time $\sigma_k$, the arrow does not have a $\delta$ at its tip with probability $(\beta_1 - \beta_2)/\beta_1 > 0$. Moreover, if $A_{\sigma_k} \notin J$ and the arrow does not have a $\delta$ at its tip, then with probability at least $1/2d > 0$, $A_s$ moves closer to 0. In particular, if $N = \inf\{k \geq 1 : A_{\sigma_k} \in J\}$, then there is $c > 0$ such that

$$P(N \geq cL) \leq C_2 \exp(-\gamma_2 L)$$

for suitable $C_2, \gamma_2 \in (0, \infty)$. Since $P(\sigma_k - \sigma_{k-1} > t) = \exp(-\beta_1 t)$, the result follows. □

LEMMA 3.3. *Assume that $\beta_1 > \beta_2$. For any $y \in J$ there exist $C_3 < \infty$ and $\gamma_3 > 0$ such that*

$$\sup_{y \in J} P_y(A_s \notin B \text{ for some } s \leq T) \leq C_3 \exp(-\gamma_3 L)$$

*for $L$ sufficiently large.*

PROOF. We let $s_0 = 0$ and, for $k \geq 1$, define the following stopping times:

$$t_k = \inf\{t > s_{k-1} : A_t \notin (-L/2, L/2)^d\},$$
$$s_k = \inf\{t > t_k : A_t \in J\},$$
$$\tau = \inf\{t > 0 : A_t \notin B\}.$$

Moreover, we denote by $M(t) = \sup\{k \geq 1 : \sigma_k < t\}$ the number of tips of arrows encountered by $A_s$ by time $t$. Then for any site $y \in J$

$$\begin{aligned}
P_y(A_s &\notin B \text{ for some } s \leq T) \\
&= P_y(A_{\sigma_k} \notin B \text{ for some } k \leq M(T)) \\
&\leq P_y(A_{\sigma_k} \notin B \text{ for some } k \leq 2\beta_1 T) + P(M(T) > 2\beta_1 T) \\
&\leq P_y(s_k > \tau \text{ for some } k \leq 2\beta_1 T) + P(M(T) > 2\beta_1 T) \\
&\leq 2\beta_1 T \sup_{z \in J} P_z(s_1 > \tau) + P(M(T) > 2\beta_1 T).
\end{aligned}$$

Since $A_s$ has a drift toward $J$ and the time between consecutive jumps has exponential bound, $P(s_1 > \tau) \leq C_4 \exp(-\gamma_4 L)$ for appropriate $C_4 < \infty$ and $\gamma_4 > 0$ (see the proof of Lemma 3.2). Furthermore, since $\mathbb{E}M(T) = \beta_1 T$, large deviation estimates imply that there are $C_5 < \infty$ and $\gamma_5 > 0$ such that $P(M(T) > 2\beta_1 T) \leq C_5 \exp(-\gamma_5 T)$. □

LEMMA 3.4. *Assume that $x \in B_z$, $\|z\| = 1$ and $\beta_1 > \beta_2$. There exist $C_6$, $\gamma_6 \in (0, \infty)$ such that*

$$\sup_{y \in B_z} P_x(A_T \notin B) \leq C_6 \exp(-\gamma_6 L)$$

*for $\Gamma$ and $L$ sufficiently large.*

PROOF. By decomposing according to whether $A_s \in J$ for some $s \leq T$ or not, we obtain

$$P_x(A_T \notin B) \leq P_x(A_s \notin J \text{ for all } s \leq T) + P(A_T \notin B; A_s \in J \text{ for some } s \leq T).$$

The first term on the right-hand side can be bounded using Lemma 3.2. For the second term, we first observe that

$$P(A_T \notin B; A_s \in J \text{ for some } s \leq T) \leq \sup_{y \in J} P_y(A_s \notin B \text{ for some } s \leq T)$$

and then apply Lemma 3.3. This completes the proof. □

Since there are $(2L+1)^d$ sites in $B_z$, it follows from Lemma 3.4 and duality that there is a constant $C_7 < \infty$ independent of $L$ such that for $\Gamma$ and $L$ sufficiently large

$$\begin{aligned}
P(\eta_T(x) = 2 \text{ for some } x \in B_z) &\leq \sum_{x \in B_z} P_x(A_T \notin B) \\
&\leq (2L+1)^d C_6 \exp(-\gamma_6 L) \leq C_7 L^{-1} \leq \varepsilon/3.
\end{aligned}$$

Moreover, since each of the dual paths has a drift toward $J$, we can fix $M > 0$, say $M = 3$, so that for any $\varepsilon > 0$

$P($any of the selected paths is not contained in

$[-ML, ML]^d$ at some time $s \leq T) \leq \varepsilon/3$

by choosing $L$ sufficiently large. This shows that boxes that are sufficiently far apart are independent of each other with high probability and completes the proof of Proposition 3.1.

**4. Proof of Theorem 1.** This section is devoted to the proof of part (b) of Theorem 1 which describes the behavior of the process for $\beta > 1$ and $g$ close to 0 in the generalist case $\alpha = \beta$. As already explained in the Introduction, if $g = 0$, then the process $\zeta_t : \mathbb{Z}^d \longrightarrow \{0, 1\}$ reduces to the biased voter model with parameters 1 and $\beta$ so that if $\beta > 1$, then $P(\zeta_t(x) = 1) \to 1$ as $t \to \infty$. To prove that the pathogens still survive when $g > 0$ is sufficiently small, we show that for $M = 3$ and any $\varepsilon > 0$, we can choose $L$ and $\Gamma$ such that Proposition 3.1 holds for $\beta_1 = \beta$ and $\beta_2 = 1$. Now that $\Gamma$ and $L$ are fixed and $M = 3$, we can assert that there exists $g_c > 0$ small so that

$P($some infected host in $[-ML, ML]^d$ gives birth

to a healthy host sent to $B$ between time 0 and time $T)$

$\leq (2ML + 1)^d (1 - \exp(-g_c T)) \leq \varepsilon/3.$

This tells us that if $g < g_c$, then the set of occupied sites dominates the set of open sites in an oriented percolation process with parameter $p = 1 - \varepsilon$. Here $(z, n)$ occupied means that all sites in $B_z$ are occupied by pathogens at time $nT$. Finally, by taking $\varepsilon > 0$ sufficiently small so that percolation occurs with positive probability, Theorem 1 follows.

**5. Proof of Theorem 2.** The proof of Theorem 2 also relies on a perturbation argument. In the case $\beta = \infty$ and $R_1 \leq R_2$, the transition $(1, 0) \to (1, 1)$ is instantaneous, that is, unassociated hosts of type 1 become instantaneously associated with a mutualist, provided that all the hosts of type 1 are initially associated with a mutualist. The assumption $R_1 \leq R_2$ is to avoid the problem of births of isolated, unassociated hosts of the same type that are not accessible to mutualists. Under these assumptions, the process $\xi_t^1$ performs a biased voter model with parameters $\beta_1 = g$ and $\beta_2 = 1$. In particular, well-known results about the biased voter model imply that if $g > 1$, then $P(\xi_t^1(x) = 1) \to 1$ as $t \to \infty$.

To extend the result to the region $\beta > 0$ large, we prove that if hosts of type 1 become occupied by their associated mutualists quickly enough, then $\xi_t^1$ will evolve like a biased voter model in the space-time box $B \times [0, \Gamma L]$ with probability close to 1. We first define $\xi_t$ on the same space as the biased voter model $\eta_t$ introduced in Section 3 with $\beta_1 = g$ and $\beta_2 = 1$. At time $T_n^{x,w}$, the host present at site $w$ gives birth to an unassociated host of the same type which is then sent to $x$. At time $U_n^{x,w}$, the birth from $w$ to $x$ occurs only if the host at $w$ is associated with a mutualist. To describe the evolution of the mutualists, we consider one more collection of independent Poisson processes, $\{V_n^{x,w} : n \geq 1\}$, $0 < \|x - w\| \leq R_2$, with parameter $\beta$. At time $V_n^{x,w}$, we draw an arrow labeled with a 1 from $w$ to $x$ to indicate that a mutualist (of type 1) present at site $w$ gives birth to a mutualist at site $x$ if this site is already occupied by a host of type 1. We will prove that there exists $\beta_{cr}^{\text{Th}2} \in (0, \infty)$ such that if $\beta > \beta_{cr}^{\text{Th}2}$ and $\xi_0^1 = \eta_0$ on $B$, then $\xi_T^1 = \eta_T$ on $B_z$ with $\|z\| = 1$ at time $T = \Gamma L$ with probability $\geq 1 - \varepsilon/3$. Since boxes that are sufficiently far apart are independent of each other with probability close to 1, we can focus on $[-ML, ML]^d \times [0, \Gamma L]$, $M = 3$, to estimate this event. Let $x \in [-ML, ML]^d$ and follow the line $\{x\} \times [0, \Gamma L]$ by going forward in time. Each time a host at $w$ attempts to give birth at site $x$, we require that the next 1-arrow from $w$ to $x$ appears before the host at $w$ is replaced or the host at $x$ gives birth. A straightforward calculation shows that this event occurs with probability

$$P(V_1^{x,w} < \min(T_1^{y,x}, U_1^{y,x}) \text{ for any } y \in \mathcal{N}_x^1 \text{ and}$$
$$V_1^{x,w} < \min(T_1^{w,y}, U_1^{w,y}) \text{ for any } y \in \mathcal{N}_w^1) = \beta(\beta + 2m)^{-1},$$

where $m = g\nu_{R_1}$. Let us now denote by $K(x,T)$ the number of unlabeled arrows and $\delta$-arrows that point at site $x$ by time $T$, and set $I_M = [-ML, ML]^d$. Then, by observing that $\mathbb{E}K(x,T) = mT$, and by decomposing the event to be estimated according to whether $K(x,T) > 2mT$ or not, we finally obtain

$$P(\xi_T^1 \neq \eta_T \text{ on } B_z) \leq \sum_{x \in I_M} P(K(x,T) > 2mT) + 2mT \sum_{x \in I_M} \frac{2m}{\beta + 2m}$$
$$\leq (2ML)^d \times \{C_8 \exp(-\gamma_8 T) + 4m^2 T(\beta + 2m)^{-1}\}$$

for appropriate $C_8 < \infty$ and $\gamma_8 > 0$. Taking $L$ and then $\beta$ sufficiently large so that

$$P(\xi_T^1 \neq \eta_T \text{ on } B_z) \leq \varepsilon/3,$$

and applying Proposition 3.1 imply that the set of occupied sites dominates the set of open sites in an oriented percolation process with parameter $p = 1 - \varepsilon$. Here $(z, n)$ occupied means that all sites in $B_z$ are occupied by associated hosts of type 1 at time $nT$. This almost produces Theorem 2. Our last problem is that oriented site percolation has a positive density of unoccupied sites. To prove that there is an in-all-directions expanding region which is void of hosts of type 2, we apply a result from Durrett [7] which shows that unoccupied sites do not percolate when $\varepsilon$ is close enough to 0. Since hosts of either type cannot appear spontaneously, once a region is void of one type, this type can only reappear in the region through invasion from the outside. This implies that $P(\xi_t^1(x) = 2) \to 0$ as $t \to \infty$ for any $x \in \mathbb{Z}^d$ and completes the proof of Theorem 2.

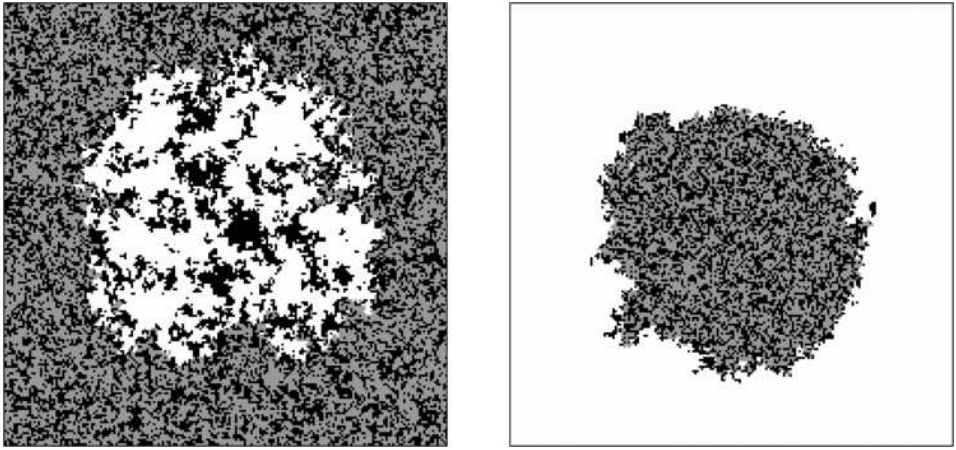

FIG. 3. *Process with nearest-neighbor interactions on the $200 \times 200$ torus at time $250$. Left $\beta = 2$ and $g = 0.5$. The process starts with unassociated white hosts in $J = (90, 110)^2$ and black hosts associated with a pathogen outside $J$. Right: $\beta = 4$ and $g = 2$. The process starts with black hosts associated with a mutualist in $J$ and unassociated white hosts outside $J$. In both pictures, gray sites refer to black hosts associated with a symbiont.*

**6. Proof of Theorem 3.** This section is devoted to the proof of Theorem 3. The method of the proof can also be applied to give a more explicit proof of Theorem 2 without too much more work. For the proof, we will assume that the birth rate $\lambda$ is not set equal to 1. In fact, we will prove Theorem 3 first for $\lambda$ small and then change the time scale so that it holds for $\lambda = 1$ as well. We start by introducing the rightmost host-2 process $r_t^2$ and the leftmost symbiont-1 process $\ell_t^1$, that is,

$$r_t^2 = \sup\{x \in \mathbb{Z} : \xi_t^1(x) = 2\} \quad \text{and} \quad \ell_t^1 = \inf\{x \in \mathbb{Z} : \xi_t^2(x) = 1\}.$$

We observe that because of nearest-neighbor interactions, $\ell_t^1 - r_t^2 \geq 1$ at any time $t \geq 0$. Moreover, if $G_t = \ell_t^1 - r_t^2 - 1$ denotes the number of sites between both processes, then $r_t^2$ changes as follows:

$$\text{for } G_t = 0 : r_t^2 \to \begin{cases} r_t^2 + 1, & \text{at rate } \lambda, \\ r_t^2 - 1, & \text{at rate } \lambda g, \end{cases}$$

and

$$\text{for } G_t > 0 : r_t^2 \to \begin{cases} r_t^2 + 1, & \text{at rate } \lambda, \\ r_t^2 - 1, & \text{at rate } \lambda, \end{cases}$$

which suggests that $r_t^2$ drifts to the right if $g < 1$ (and to the left if $g > 1$). To make this argument precise, we will prove the following lemma.

LEMMA 6.1. *For $g \neq 1$, there exists $\varepsilon_0 > 0$ such that $\lim_{T \to \infty} \frac{1}{T} \times \int_0^T \mathbb{1}_{\{G_t = 0\}} \, dt \geq \varepsilon_0$.*

To deduce Theorem 3 from Lemma 6.1, we first observe that $\mathbb{E} r_t^2 \geq \lambda(1 - g)\varepsilon_0 t$ for $g < 1$. A large deviation estimate then implies that

$$P(r_t^2 \geq (1-\varepsilon)\sigma t) \geq 1 - C_9 e^{-\gamma_9 t},$$

where $\sigma = \lambda(1-g)\varepsilon_0$. This implies that if $g < 1$, then $r_t^2 \to \infty$ with probability 1. If $g > 1$, then $\mathbb{E} r_t^2 \leq -\lambda(g-1)\varepsilon_0 t$, and a similar argument implies that $r_t^2 \to -\infty$ with probability 1. The proof of Lemma 6.1 is based on a rescaling argument. The main objective is to prove that the fraction of time the host present at site $r_t^2 + 1$ is associated with a symbiont is greater than some positive constant. To be able to compare the particle system with oriented percolation process, we will artificially freeze $r_t^2$ by introducing the process seen from the interface, namely $\iota_t(x) = \xi_t(x + r_t^2 + 1)$. From this new point of view, the symbiont evolves on the half-line $\mathbb{Z}^+$. To do this comparison, we need to extend some results on oriented percolation to oriented percolation *in half-space*, that is, the process in which all sites to the left of 0 are closed.

*Oriented percolation in half-space.* As explained above, our first objective is to investigate 1-dependent site percolation process in the half-space. Let

$$\Omega = \{(x,n) : x + n \text{ is even and } n \geq 0\}.$$

For any $(x, n) \in \Omega$, let $\omega(x, n)$ define a random variable with the following property: If $x < 0$, then $\omega(x, n) = 0$, while if $x \geq 0$, then $\omega(x, n) = 1$ with probability $p$ and $\omega(x, n) = 0$ with probability $1 - p$. The site $(x, n)$ is said to be open (resp. closed) if $\omega(x, n) = 1$ (resp. 0). Finally, 1-dependent with parameter $p$ means that whenever $(x_i, n_i)$, $1 \leq i \leq m$, is a sequence with $x_i \geq 0$ for any $1 \leq i \leq m$ where $\|(x_i, n_i) - (x_j, n_j)\|_\infty > 1$ for $i \neq j$, then

$$P(\omega(x_i, n_i) = 0 \text{ for } 1 \leq i \leq m) = (1-p)^m.$$

We say that $(y, n)$ can be reached from $(x, m)$ and write $(x, m) \to (y, n)$ if there exists a sequence of points $x = x_m, x_{m+1}, \ldots, x_n = y$ such that $|x_{i+1} - x_i| = 1$ for $m \leq i \leq n-1$, and $\omega(x_i, i) = 1$ for $m \leq i \leq n$. To formulate the next result, we let

$$W_n^A = \{y : (x, 0) \to (y, n) \text{ for some } x \in A\} \quad \text{and} \quad \tau^A = \inf\{n : W_n^A = \varnothing\}$$

for any $A \subset 2\mathbb{Z}^+$. Finally, if $A = 2\mathbb{Z}^+$, we will write $W_n = W_n^A$ and $\tau = \tau^A$.

LEMMA 6.2. *If $p > 1 - 3^{-72}$, then there exists $\varepsilon_1 > 0$ such that $P(0 \in W_{2n}) \geq \varepsilon_1$ for any $n \geq 0$.*

PROOF. The proof closely follows ([5], Section 10) so we will try to be as brief as possible, and refer the reader to this reference for details. The first step is a slight modification of the contour argument applied to independent site percolation process. Let $A = \{0, 2, \ldots, 2N\}$ and $C = \{(y, n) : \text{there exists } x \in A \text{ with } (x, 0) \to (y, n)\}$. We set

$$D = \{(a, b) \in \mathbb{R}^2 : |a| + |b| \leq 1\} \quad \text{and} \quad W = \bigcup_{z \in C}(z + D).$$

If the set $C$ is finite, we denote by $\Gamma$ the boundary of the unbounded component of $(\mathbb{R} \times (-1, \infty)) - W$ and orient the boundary in such a way that the segment from $(2N, -1)$ to $(2N + 1, 0)$ is oriented in the direction indicated. The boundary is a contour line, if it exists, starting at $(2N, -1)$ and ending at $(0, -1)$. There are at most $3^{m-1}$ contours of length $m$. Moreover, for a contour of length $m$ to exist, there must be at least $m/8$ closed sites to the right of zero. To prove this point, we call a segment a line segment of the form $x + F$ where $x \in C$ and $F$ is one of the sides of $D$. The site closest to the right of the segment is the site associated with the segment. Here, right and left are defined according to the orientation introduced above. We call segments of $\Gamma$ which look like $\nwarrow$, $\nearrow$, $\searrow$ and $\nearrow$ segments of types 1, 2, 3 and 4, respectively. By construction, a site associated with a segment of type 1 or 2 must be closed. Let $m_i$ and $\bar{m}_i$ be, respectively, the number of segments of type $i$ and the number of segments of type $i$ located on the left of zero. Since the contour starts at $(2N, -1)$ and ends at $(0, -1)$, $m_1 + m_2 = m_3 + m_4 + 2N$, so if the contour has length $m$, then $m_1 + m_2 \geq m/2$. The same reasoning leads to $m_1 + m_4 = m_2 + m_3$ so that $\bar{m}_2 + \bar{m}_3 \leq m/2$. Now, since the sites located on the left of zero are closed with probability 1, we obtain $\bar{m}_1 = 0$ and $\bar{m}_2 = \bar{m}_3 \leq m/4$, which implies $(m_1 + m_2) - (\bar{m}_1 + \bar{m}_2) \geq m/4$. Finally, since a site in $W^c$ can be associated with at most two segments of type 1 and 2, it follows that the number of sites on the right of zero that must be closed is $\geq m/8$. Noticing that the shortest possible contour has length $2N + 4$, one can conclude that

$$P(\tau^{[0, 2N]} < \infty) \leq \sum_{m=2N+4}^{\infty} 3^m (1-p)^{m/8} = C_{10}(3(1-p)^{1/8})^{2N}$$

if $p > 1 - 3^{-8}$ and the variables $\omega(x, n)$ are independent. To deal with the 1-dependent case, we observe that there are nine sites in $\Omega$ with $(|m| + |n|)/2 \leq 1$, so for each $\Gamma$ of length $m$ there is a set of $m/72$ sites which are separated by more than one and which must be closed for the contour to exist. In conclusion, if $p > 1 - 3^{-72}$, then $P(\tau^0 = \infty) > 0$. Finally, if we map

$\Omega$ into itself by sending $(x, m) \mapsto (x, n - m)$ and reverse the orientation of the graph, then

$$P(W_n^A \cap B \neq \varnothing) = P(W_n^B \cap A \neq \varnothing).$$

Taking $A = 2\mathbb{Z}^+$ and $B = \{0\}$, it then follows that $P(0 \in W_{2n}) \downarrow P(\tau^0 = \infty) > 0$. In conclusion, the lemma holds by setting $\varepsilon_1 = P(\tau^0 = \infty)$. □

*The process seen from the interface.* To prove Lemma 6.1 for $g < 1$, we introduce the nearest-neighbor contact process $\xi_t^-$ in which a particle at $x$ dies at rate $2\lambda$, gives birth at rate $\beta \nu_{R_2}$ and sends its offspring to one of the neighbors at random in $\mathcal{N}_x^2$. The process is modified so that particles located in the interval $(-\infty, r_t^2]$ are removed. More precisely, each time a particle in $\xi_t^-$ tries to give birth to a particle which is sent to a site in $(-\infty, r_t^2]$, the birth is suppressed. Furthermore, if a particle is present at site $r_t^2 + 1$ when the process $r_t^2$ moves to the right, then this particle is removed. (The proof is similar in the case $g > 1$, provided one replaces $\xi^-$ by $\xi^+$ where particles give birth at rate $\beta \nu_{R_2}$ but now die at rate $2\lambda g$.) The first step is to prove that $\xi_t$ dominates the process $\xi_t^-$.

LEMMA 6.3. *If $g < 1$, the processes $\xi_t$ and $\xi_t^-$ can be defined on the same probability space in such a way that if the inclusion $\xi_0^- \subset \xi_0^2$ holds at time 0, then $\xi_t^- \subset \xi_t^2$ at any later time.*

PROOF. Let us start by observing that if the rightmost process $r_t^2$ jumps to the right, a particle located at $r_t^2 + 1$, if it is present, is removed from both processes $\xi_t^-$ and $\xi_t^2$ and that, if it jumps to the left, both configurations stay unchanged. Therefore, it suffices to prove that the inclusion holds when $r_t^2$ is constant, say $r_t^2 = 0$ at any time. This follows from a standard coupling argument so we just need to define a graphical representation that preserves the inclusion. We consider for any $x, z \geq 0$ with $|x - z| = 1$ the independent Poisson processes $\{S_n^{x,z} : n \geq 1\}$ and $\{T_n^{x,z} : n \geq 1\}$ with rate $\lambda g$ and $\lambda(1-g)$, respectively. For any $x, z \geq 0$, with $0 < \|z - x\| \leq R_2$, we also introduce the Poisson process $\{U_n^{x,z} : n \geq 1\}$ with rate $\beta$. The evolution of $\xi_t$ is as follows: At time $S_n^{x,z}$, the host present at site $x$ gives birth to a healthy host of the same type, which is then sent to $z$. At time $T_n^{x,z}$, the birth occurs only if the host at $x$ is healthy. At time $U_n^{x,z}$, a pathogen (of type 1) at site $x$ infects a host of type 1 at $z$ if it is present. Finally, the process $\xi_t^-$ evolves according to the following rules: At time $U_n^{x,z}$, a particle at site $x$ gives birth to a new particle which is then sent to $z$. If the site is empty, the birth occurs. Otherwise, it is suppressed. At times $T_n^{x,z}$ and $S_n^{x,z}$, a particle present at $z$ is removed. Such a coupling leads to the desired result. □

We now introduce the process seen from the interface: $\eta_t(x) = \xi_t^-(x + r_t^2 + 1)$. To describe this process, we define the translation operators $\tau_1$ and $\tau_{-1}$ by setting

$$[\tau_1 \eta](x) = \eta(x+1) \quad \text{and} \quad [\tau_{-1}\eta](x) = \begin{cases} \eta(x-1), & \text{if } x \geq 1, \\ 0, & \text{if } x \leq 0. \end{cases}$$

Then $\eta_t$ evolves according to the following rules:

1. A particle at $x$ gives birth at rate $\beta \nu_{R_2}$ to a new particle which is then sent to a neighbor within the neighborhood $\mathcal{N}_x^2$. If the target site is empty the birth occurs; otherwise it is suppressed. Moreover, no births are allowed to the left of 0.
2. Each particle dies at rate $2\lambda$.
3. Depending on whether $\xi_t^2(r_t^2 + 1) = 0$ or $1$, $\eta_t$, respectively, shifts as follows:

$$\eta_t \to \begin{cases} \tau_1 \eta_t, & \text{at rate } \lambda, \\ \tau_{-1}\eta_t, & \text{at rate } \lambda, \end{cases} \quad \text{and} \quad \eta_t \to \begin{cases} \tau_1 \eta_t, & \text{at rate } \lambda, \\ \tau_{-1}\eta_t, & \text{at rate } \lambda g. \end{cases}$$

To compare $\eta_t$ with an oriented percolation process in half-space, we let $\Gamma$ and $L$ be two large integers to be fixed later, and scale space by setting

$$B = [-L, L], \qquad \Phi(z) = Lz \quad \text{and} \quad B_z = \Phi(z) + B.$$

Let $J_z = \Phi(z) + (-L/5, L/5)$, and say that the site $(z, n) \in \Omega$ is *occupied* if there is at least one particle in the interval $J_z$ at time $n\Gamma L$. Let us fix $z \geq 0$ even, and start the process $\eta_t$ with one particle in $J_z$. The first step is to investigate the process with no deaths inside $B_z$ and modified so that any particle outside $B_z$ is killed. We denote by $\bar{\eta}_t$ this new process.

LEMMA 6.4. *Let $\varepsilon_2 = 6^{-72}$ and $\beta > \lambda$. Then $\Gamma$ and $L$ can be chosen so that*

$$P(\bar{\eta}_{\Gamma L} \cap J_{z+1} \neq \varnothing) \geq 1 - 2\varepsilon_2/3.$$

PROOF. A standard coupling argument implies that $\bar{\eta}_t$ has more particles if we increase the range of the interactions. So, we just need to prove the result when the offspring is sent to one of the two nearest neighbors. The idea is to prove that the rightmost particle in $\bar{\eta}_t$ reaches the right edge of $B_z$ by time $\Gamma L$, and then stays inside $J_{z+1} \cap B_z$ until time $\Gamma L$. Let

$$r_t = \sup\{x \in B_z : \bar{\eta}_t(x) = 1\} \quad \text{and} \quad \tau = \inf\{t \geq 0 : r_t = \Phi(z+1)\}.$$

Then, on the set $\{\tau > t\}$, we have $\mathbb{E}(r_t - r_0) \geq (\beta - \lambda)t$. In particular, since $\beta > \lambda$, the parameters $\Gamma$ and $L$ can be chosen such that $P(\tau > \Gamma L) \leq \varepsilon_2/3$. This implies that the rightmost particle will reach $J_{z+1}$ by time $\Gamma L$ with high probability. To prove that the rightmost particle does not leave $J_{z+1}$ until time $\Gamma L$, we observe that

$$r_t \to r_t - 1 \quad \text{at rate} \quad \leq \lambda \quad \text{and}$$

$$r_t \to r_t + 1 \quad \text{at rate} \quad \begin{cases} \geq \beta, & \text{if } r_t < \Phi(z+1), \\ 0, & \text{if } r_t = \Phi(z+1). \end{cases}$$

Then well-known estimates about random walks imply that

$$P(\exists t \in [\tau, \Gamma L] : r_t \notin J_{z+1}) \leq \Gamma L \left(\frac{\lambda}{\beta}\right)^{L/5} \leq \varepsilon_2/3$$

for $L$ sufficiently large. This completes the proof. $\square$

We now fix $\Gamma$ and $L$ such that Lemma 6.4 holds. To extend the result to the process $\eta_t$, we just need to choose $\lambda > 0$ sufficiently small so that the probability a death occurs in the space-time region $B_z \times [0, \Gamma L]$ is smaller than $\varepsilon_2/3$. In other respects, since the result holds for the process modified so that any particle outside $B_z$ is killed, it follows that: If $\lambda > 0$ is small, then the set of occupied sites dominates the set of wet sites in an oriented percolation process in the half-space with parameter $1 - 6^{-72}$. Lemma 6.2 then implies that

$$P(\text{there is at least one particle in } J_z \text{ at time } 2n\Gamma L) \geq \varepsilon_1 > 0$$

for any integer $n \geq 0$, provided that $\eta_0$ contains infinitely many particles. Now, it is easy to see that there exists a constant $\varepsilon_3 > 0$ independent of $n$ such that: If there is at least one particle in the interval $J_z$ at time $2n\Gamma L$, then the probability that $\xi_t^-(r_t^2 + 1) = \eta_t(0) = 1$ for at least one unit of time between times $2n\Gamma L$ and $2(n+1)\Gamma L$ is greater than $\varepsilon_3$. This tells us that

$$\lim_{T \to \infty} \frac{1}{T} \int_0^T \mathbb{1}_{\{\eta_t(0)=1\}} \, dt \geq \varepsilon_0$$

for some appropriate constant $\varepsilon_0 > 0$. Since Lemma 6.3 implies that $\{\eta_t(0) = 1\} \subset \{G_t = 0\}$, Lemma 6.1 follows from the previous inequality. A time change now allows us to set $\lambda = 1$, which then completes the proof of Theorem 3.

**7. Proof of Theorem 4.** This section is devoted to the proof of Theorem 4 which addresses coexistence of the symbionts in the neutral case $g = 1$. To remind the reader, we assume that the symbionts evolve as previously but the hosts perform a threshold $\theta$ voter model according to the following rate at $x$:

$$i \to j \quad \text{at rate} \quad \begin{cases} 1, & \text{if } \text{card}\{z \in \mathbb{Z}^d : 0 < \|x - z\| \leq R_1 \\ & \quad \text{and } \hat{\xi}_t^1(z) = j\} \geq \theta, \\ 0, & \text{otherwise.} \end{cases}$$

It is easy to see that the critical value for the infection rate, $\beta_{cr}^{\text{Th}4}(\kappa)$, is strictly bounded away from 0. Namely, if $\kappa = 1$, then the symbionts perform a basic contact process with death rate 1, provided $\theta \leq \nu_{R_1}$, and birth rate $\beta \nu_{R_2}$. Furthermore, since the contact process is monotone, $\beta_{cr}^{\text{Th}4}(\kappa) \geq \beta_{cr}^{\text{Th}4}(1)$ for $\kappa \geq 1$, from which our claim follows.

To prove Theorem 4, we will compare the particle system viewed on suitable length and time scales with a 1-dependent oriented percolation process in two dimensions. The properties of the process in the absence of symbionts was described in [7]. To apply his results, we introduce, for any $x \in \mathbb{Z}^d$, the *house*

$$H_x = [x_1 L, (x_1 + 1)L) \times \cdots \times [x_d L, (x_d + 1)L),$$

where $L$ is an integer to be fixed later and $x_i$ denotes the $i$th coordinate of the vector $x$. We fix $\sigma < 1/\kappa$ such that $\theta < \sigma \nu_{R_1}$, and say that $H_x$ is *good* if it contains at least $\sigma L^d$ hosts of each type. For $x \in \mathbb{Z}^d$, we define $\|x\|_2 = (|x_1|^2 + \cdots + |x_d|^2)^{1/2}$ and set $B_2(x, r) = \{y : \|y - x\|_2 \leq r\}$. We say that $B_2(0, r)$ is *good* if for any $x \in B_2(0, r)$ the house $H_x$ is good. For $z$ even for even $n$ or $z$ odd for odd $n$, we will say that $(z, n)$ is *occupied* if the following two conditions hold:

1. For any $x \in B_2(zKe_1, K)$, the house $H_x$ is good at time $n\Gamma L$.
2. For any $i = 1, 2, \ldots, \kappa$, $B_2(zKLe_1, KL)$ has at least one symbiont of type $i$ at time $n\Gamma L$.

Here, $e_1$ denotes the first unit vector, and $K$ and $\Gamma$ are large integers that will be fixed later. Note that the set $B_2(zKe_1, K)$ is defined on the rescaled lattice, whereas $B_2(zKLe_1, KL)$ is defined on the original lattice. We will prove the following result.

PROPOSITION 7.1. *Let $\varepsilon > 0$ and $\theta < \nu_{R_1}/\kappa$. There exists $\beta_{cr}^{\text{Th}4} \in (0, \infty)$ such that if $\beta > \beta_{cr}^{\text{Th}4}$, then $K$, $L$ and $\Gamma$ can be chosen in such a way that the set of occupied sites dominates the set of open sites in a 1-dependent oriented percolation process with parameter $1 - \varepsilon$.*

The first step in proving Proposition 7.1 is to summarize the results of Durrett ([7], Section 2), which describe the behavior of the process in the absence of symbionts. To formulate the result we are interested in, we set $R_1 = L(M_1 + 1)$ where $L$ and $M_1$ are large integers.

LEMMA 7.2 (Durrett). *Let $\varepsilon > 0$ and $\theta < \nu_{R_1}/\kappa$. There exist $R_0$, $M_0$ and $\Gamma$ such that the following holds: If $M_1 \geq M_0$ and $B_2(0, R_0 M_1)$ is good at time 0, then, for $L$ large, $B_2(0, R_0 M_1)$ is good until time $\Gamma L$ and $B_2(0, 2R_0 M_1)$ is good at time $\Gamma L$ with probability at least $1 - \varepsilon/3$.*

The sets described in Lemma 7.2 provide an environment favorable to the survival of symbionts. To explain this, we introduce, for any type $i = 1, 2, \ldots, \kappa$, the processes $\eta_t^i$ defined by $\eta_t^i(x) = 1$ if $\hat{\xi}_t^2(x) = i$ and $\eta_t^i(x) = 0$ otherwise. Since $\alpha = 0$, it is easy to see that, for $i = 1, 2, \ldots, \kappa$, the processes $\eta_t^i$ do not interact. We fix a type $i \in \{1, 2, \ldots, \kappa\}$, and focus on the process $\eta_t^i$. The evolution of $\eta_t^i$ is as follows:

1. Each particle dies at rate at most $\kappa$ and gives birth at rate $\beta \nu_{R_2}$.
2. A particle born at site $x$ is sent to a site $z$ chosen at random from $\mathcal{N}_x^2$.
3. If the target site $z$ is occupied by an unassociated host of type $i$, then the birth occurs. Otherwise, it is suppressed.

The proof of Theorem 4 relies, like Theorems 1 and 2, on a perturbation argument. More precisely, we first prove Proposition 7.1 in the extreme case $\beta = \infty$, and then extend the result to the region $\beta > 0$ large. We denote by $\bar{\eta}_t^i$ the process $\eta_t^i$ modified so that no births are allowed outside $B_2(0, KL)$.

LEMMA 7.3. *Assume that $B_2(0, K)$ is good until time $\Gamma L$ and that at time 0 there exists $x \in B_2(0, KL)$ with $\eta_0^i(x) = 1$. If $R_2 \geq 4dL$ and $\beta = \infty$, then*

$$\{x \in \mathbb{Z}^d : \hat{\xi}_t^1(x) = i\} \cap B_2(0, KL) = \{x \in \mathbb{Z}^d : \bar{\eta}_t^i(x) = 1\} \qquad \text{for all } t \leq \Gamma L.$$

PROOF. This is elementary geometry. To begin with, we cover the set $B_2(0, KL)$ with a finite number of Euclidean balls $B_j$, $j \in I$; each of them has radius $r = \sqrt{d}L$. Then, it is easy to see that, for any $j \in I$, $B_j$ contains at least one house. In particular, as long as $B_2(0, K)$ is good, $B_j$ contains at least one host of type $i$ provided that $B_j \subset B_2(0, KL)$. At any time $0 \leq t \leq \Gamma L$, let us pick one at random and denote by $X_j(t)$ its spatial location. Now, since $R_2 \geq 4r\sqrt{d}$, we have

$$\min\{\|X_j(t) - X_k(t)\| : k \neq j\} \leq R_2 \qquad \text{for all } j \in I \text{ and for all } t \leq \Gamma L.$$

This implies that, for any $x, z \in B_2(0, KL)$ occupied by a host of type $i$, there exists a chain of sites $x_0 = x, x_1, \ldots, x_n = z$ such that the following two conditions hold:

1. For $k = 1, 2, \ldots, n$, $\|x_{k-1} - x_k\| \leq R_2$.
2. For $k = 0, 1, \ldots, n$, the site $x_k$ is occupied by a host of type $i$.

In particular, since $\bar{\eta}_t^i$ starts with at least one particle in $B_2(0, KL)$ and $B_2(0, KL)$ is finite, all the hosts of type $i$ are instantaneously invaded by a symbiont at time 0. It is easy to prove by induction that this holds until time $\Gamma L$. If a host of type $i$ gives birth to an unassociated host which is sent to a site $x \in B_2(0, KL)$ at time $t$, pick $X_j(t)$ such that $\|x - X_j(t)\| \leq$ $R_2$. Since $X_j(t)$ is occupied by a symbiont of type $i$, the host at $x$ will be instantaneously invaded. □

To extend the result to $\beta > 0$ large, it is convenient to construct the process $\hat{\xi}_t$ from a graphical representation. For any type $i \in \{1, 2, \ldots, \kappa\}$ and $x \in \mathbb{Z}^d$, let $\{T_n^{i,x} : n \geq 1\}$ be independent Poisson processes with rate 1. At time $T_n^{i,x}$ the state of $x$ flips to $(i, 0)$ if the set $\mathcal{N}_x^1$ has at least $\theta$ hosts of type $i$. For $x \in \mathbb{Z}^d$, let $\{U_n^x : n \geq 1\}$ be independent Poisson processes with rate $\beta$. At time $U_n^x$, we choose at random a site $z$ from $\mathcal{N}_x^2$. If a host of a certain type is present at site $x$, and a symbiont of the same type is present at site $z$, then the host at site $x$ becomes associated if it is not already. So that Lemma 7.3 holds for $\beta < \infty$ large, we now require the following two good events, denoted by $G_1$ and $G_2$, respectively: First, we need a *quick invasion* of the ball $B_2(0, KL)$ by the symbionts. More precisely, $G_1$ will be the event that if at time 0 there exists $x \in B_2(0, KL)$ with $\eta_0^i(x) = 1$, then for all $z \in B_2(0, KL)$, with $z$ occupied by host $i$, the host present at site $z$ becomes associated before another host attempts to give birth in $B_2(0, KL)$. To estimate $P(G_1)$, we observe that, for any $z \in B_2(0, KL)$, the host at site $z$ can be reached in at most $4KL/R_2$ steps by a symbiont, that is, if the host at $z$ is of type $i$, then there is a chain of sites $x_0, x_1, \ldots, x_n = z$ with $n \leq 4KL/R_2$, satisfying the conditions 1 and 2 above and such that $x_0$ is occupied by a symbiont of type $i$ at time 0. We denote by $\nu_{KL}$ the number of sites in the ball $B_2(0, KL)$. Then since the transition $i \to j$ occurs at rate at most 1 and there are $\kappa$ hosts and $\nu_{KL}$ sites in $B_2(0, KL)$, new hosts are born at rate at most $\kappa \nu_{KL}$. We set $n$ equal to the integer part of $4KL/R_2$. Then

$$P(G_1) \geq 1 - \nu_{KL} n \frac{\kappa \nu_{KL}}{\kappa \nu_{KL} + \beta/n}$$

$$\geq 1 - \nu_{KL} \frac{4KL}{R_2} \frac{\kappa \nu_{KL}}{\kappa \nu_{KL} + R_2 \beta/4KL}.$$

Now that $B_2(0, KL)$ has been invaded, we secondly require it to remain *fully occupied* until time $\Gamma L$. In other words, $G_2$ will be the event that given that at time 0 all hosts are associated, each time a host is born at some site $x \in B_2(0, KL)$, it becomes associated before another host is born in the ball $B_2(0, KL)$; this occurs from time 0 to time $\Gamma L$. Let $N$ denote the number of times a host is born in $B_2(0, KL)$ from time 0 to time $\Gamma L$. Since $\mathbb{E} N \leq \kappa \nu_{KL} \Gamma L$, we find that

$$P(N > 2\kappa \nu_{KL} \Gamma L) \leq C_{11} \exp(-\gamma_{11} \Gamma L)$$

for appropriate $C_{11} < \infty$ and $\gamma_{11} > 0$. If only one host in $B_2(0, KL)$ is unassociated, it becomes associated at rate at least $\beta$. Births of hosts in $B_2(0, KL)$ occur at rate at most $\kappa\nu_{KL}$. Let $X$ be a random variable with exponential distribution with parameter $\beta$ and let $Y$ be a random variable with exponential distribution with parameter $\kappa\nu_{KL}$. Then

$$\begin{aligned}P(G_2^c) &\leq P(N > 2\kappa\nu_{KL}\Gamma L) + P(G_2^c; N \leq 2\kappa\nu_{KL}\Gamma L) \\ &\leq C_{11}\exp(-\gamma_{11}\Gamma L) + 2\kappa\nu_{KL}\Gamma L P(Y \leq X) \\ &\leq C_{11}\exp(-\gamma_{11}\Gamma L) + 2\kappa\nu_{KL}\Gamma L \frac{\kappa\nu_{KL}}{\beta + \kappa\nu_{KL}}.\end{aligned}$$

The proof of Proposition 7.1 is now straightforward. Let $\varepsilon > 0$ and assume that $B_2(0,K)$ is good and that $B_2(0,KL)$ has at least one symbiont of each type at time 0. Fix $R_1 = L(M_1 + 1)$ and $K = R_0 M_1$, then apply Lemma 7.2 and choose $L$ sufficiently large so that $B_2(0,K)$ is good from time 0 to time $\Gamma L$ and $B_2(0,2K)$ is good at time $\Gamma L$ with probability at least $1 - \varepsilon/3$. Now, increase $L$ and then choose $\beta$ sufficiently large so that both probabilities $P(G_1)$ and $P(G_2|G_1)$ are greater than $1 - \varepsilon/3$. To see that this produces the desired result, we observe that if $B_2(0,K)$ is good from times 0 to $\Gamma L$, then Lemma 7.3 implies that, on $G_1 \cap G_2$, the balls $B_2(-KLe_1, KL)$ and $B_2(KLe_1, KL)$ contain at least one symbiont of each type. This completes the proof.

To deduce the existence of a nontrivial stationary measure $\mu$ from Proposition 7.1, we start the process $\hat{\xi}_t$ from a product measure in which each host is associated with a symbiont and has density $1/\kappa$. Then, we take the Cesaro average of the distributions from time 0 to time $T$ and extract a convergent subsequence. By Proposition 1.8 of [15], the limit $\mu$ is known to be an invariant measure. To see that $\mu$ has the desired property, we observe that if $L$ is large, then the law of large numbers implies that $(z,0)$, $z$ even, is occupied with probability close to 1. Moreover, if $\varepsilon > 0$ is small, well-known percolation results imply that, at any level $n$, the density of occupied sites is positive, which implies that $\mu(\hat{\xi}^2(x) = i) \neq 0$ for any $i \in \{1, 2, \ldots, \kappa\}$. At this point, we have proved that there is a critical value $\beta_{cr}^{\text{Th}4} \in (0, \infty)$ such that if $\beta > \beta_{cr}^{\text{Th}4}$, then coexistence occurs.

To see that $\beta_{cr}^{\text{Th}4}$ can be chosen so that if $\beta < \beta_{cr}^{\text{Th}4}$ then coexistence does not occur, we rely on a standard coupling argument. If we think of the process as being generated by the Poisson processes introduced above, it is easy to see that if $\beta_1 < \beta_2$, then the processes with parameters $\beta_1$ and $\beta_2$ can be defined on the same space, starting from the same initial configuration, in such a way that the process with parameter $\beta_1$ has fewer symbionts of type $i$ for any $i \in \{1, 2, \ldots, \kappa\}$. This completes the proof of Theorem 4.

**Acknowledgment.** The authors would like to thank an anonymous referee for his/her careful reading of the proofs.

LABORATOIRE DE MATHÉMATIQUES RAPHAËL SALEM
UMR 6085, CNRS–UNIVERSITÉ DE ROUEN
AVENUE DE L'UNIVERSITÉ
BP. 12
76801 SAINT ETIENNE DU ROUVRAY
FRANCE
E-MAIL: nicolas.lanchier@univ-rouen.fr

ECOLOGY, EVOLUTION AND BEHAVIOR
UNIVERSITY OF MINNESOTA
1987 UPPER BUFORD CIRCLE
ST PAUL, MINNESOTA 55108
USA